\documentstyle[11pt]{article}
\newcommand{\qb}[2]{{\left [{#1 \atop #2} \right]}}
\newlength{\standardunitlength}
\newtheorem{cor}{Corollary} 
\newtheorem{theorem}{Theorem} \newtheorem{prop}{Proposition}
\newenvironment{proof}{\noindent {\sc Proof:}}{$\Box$ \vspace{2 ex}}
\begin{document}

\newpage

\begin{center}
The Combinatorics of Biased Riffle Shuffles
\end{center}

\begin{center}
By Jason Fulman
\end{center}

\begin{center}
Dartmouth College
\end{center}

\begin{center}
Jason.E.Fulman@Dartmouth.Edu
\end{center}

\begin{abstract}
	This paper studies biased riffle shuffles, first defined by
Diaconis, Fill, and Pitman. These shuffles generalize the
well-studied Gilbert-Shannon-Reeds shuffle and convolve nicely. An upper bound is
given for the time for these shuffles to converge to the uniform
distribution; this matches lower bounds of Lalley. A careful version of a
bijection of Gessel leads to a generating function for cycle structure after one
of these shuffles and gives new results about descents in random permutations.
Results are also obtained about the inversion and descent structure of a
permutation after one of these shuffles.
\end{abstract}

\section{Introduction and Background} \label{intro}

	The most widely used method of shuffling cards is riffle
shuffling. Roughly speaking, one cuts the deck of cards into two piles of
approximately equal size and then riffles the two piles together. A precise
mathematical model of riffle shuffles is the Gilbert-Shannon-Reeds (or GSR)
shuffle, found independently by Gilbert $\cite{Gilbert}$ and Reeds
$\cite{Reeds}$. This model says to first cut the $n$ card deck into two
packs of size $m$ and $n-m$ with probability $\frac{{n \choose
m}}{2^n}$. Then drop cards from these packs one at a time, such that if
pack 1 has $A_1$ cards and pack 2 has $A_2$ cards, the next card is dropped
from pack 1 with probability $\frac{A_1}{A_1+A_2}$ and from pack 2 with
probability $\frac{A_2}{A_1+A_2}$.

	Before defining biased shuffles, let us recall the notion of
the descent set of a permutation. An element $\pi \in S_n$ is said to
have a descent at position $i$ if $\pi(i)>\pi(i+1)$. By convention we
say that all $\pi \in S_n$ have a descent at position $n$. The descent
set of $\pi$ is the set of positions at which $\pi$ has a descent.

	This paper analyzes a notion of biased riffle shuffles which
generalizes the GSR shuffle (the GSR shuffle will correspond to the case
$a=2,p_1=p_2=\frac{1}{2}$). These biased shuffles seem to have first been
considered on pages 153-4 of Diaconis, Fill, and Pitman $\cite{DiFiPi}$. We
now give four descriptions of these biased riffle shuffles. These
descriptions generalize the descriptions of the GSR shuffle in Bayer and
Diaconis $\cite{Bayer}$. It is elementary to prove that these descriptions
are equivalent.

\begin{center}
Descriptions of Biased $a$-shuffles
\end{center}

\begin{enumerate}

\item Cut the $n$ card deck into $a$ piles by picking pile sizes
according to the $mult(a;\vec{p})$ law, where $p=(p_1,\cdots,p_a)$. In
other words, choose $b_1,\cdots,b_a$ with probability:

\[ {n \choose b_1 \cdots b_a} \prod_{i=1}^a p_i^{b_i} \] 

	Then choose uniformly one of the ${n \choose b_1 \cdots b_a}$ ways
of interleaving these packets, leaving the cards in each packet in their
original relative order. (In the language of descents, choose uniformly one
of the ${n \choose b_1 \cdots b_a}$ permutations whose inverse has descent
set contained in $\{b_1,b_1+b_2,\cdots,b_1+\cdots+b_a=n\}$).

\item As in Description 1, cut the $n$ card deck into $a$ piles
according to the $mult(a;\vec{p})$ law. Now drop cards from the $a$
packets one at a time, according to the rule that if the $i$th packet
has $A_i$ cards, then the next card is dropped from the $i$th packet
with probability $\frac{A_i}{A_1+\cdots+A_a}$.

\item Drop $n$ points in $[0,1]$ according to the following
procedure. Break the unit interval into $a$ sub-intervals of length
$\frac{1}{a}$. Pick the $i$th interval with probability proportional
to $p_i$. Then drop uniformly in this interval. Label the points
$x_1,\cdots,x_n$ in order of smallest to largest. The map $x \mapsto
ax$ (mod 1) reorders these points. The induced measure on $S_n$ is the
same as in Descriptions 1 and 2.

\item The inverse of a biased $a$-shuffle has the following
description. Start with an ordered deck of $n$ cards face
down. Successively and independently, cards are turned face up and
dealt into a random pile $i$ with probability proportional to
$p_i$. After all the cards have been distributed, the piles are
assembled from left to right and the deck is turned face down.

\end{enumerate}

	We denote the measure on $S_n$ defined by Descriptions 1-4 by
$P_{n,a;\vec{p}}$. For example, one can check that the measure
$P_{3,2;p_1,1-p_1}$ assigns to permutations in cycle form the following
masses:

\begin{eqnarray*}
(1)(2)(3) & & \ \ \ \ p_1^3+p_1^2p_2+p_1p_2^2+p_2^3\\
(1)(23) & & \ \ \ \ p_1^2p_2\\
(2)(13) & & \ \ \ \ 0\\
(3)(12) & & \ \ \ \ p_1p_2^2\\
(123) & & \ \ \ \ p_1p_2^2\\
(132) & & \ \ \ \ p_1^2p_2
\end{eqnarray*}

	If $\vec{p}=(p_1,\cdots,p_a)$ and $\vec{p'}=(p_1',\cdots,p_b')$,
define the product:

\[ \vec{p} \otimes \vec{p'}=
(p_1p_1',\cdots,p_1p_b',\cdots,p_ap_1',\cdots,p_ap_b') \]

	The following fact, which shows that biased riffle shuffles
convolve well, is stated without proof in Diaconis, Fill, and Pitman
$\cite{DiFiPi}$.

\begin{prop} \label{convol} The convolution of $P_{n,a;\vec{p}}$ and
$P_{n,b;\vec{p'}}$ is $P_{n,ab;\vec{p} \otimes \vec{p'}}$.
\end{prop}

\begin{proof}
	This follows from the inverse description of card
shuffling. Lexicographically combining the pile assignments from an inverse
a-shuffle and an inverse b-shuffles gives uniform and independent pile
assignments for an inverse ab-shuffle.
\end{proof}

	Proposition $\ref{convol}$ is the starting point for this
paper. Little seems to be known about biased riffle shuffles. The
Gilbert-Shannon-Reeds shuffle (the case of equal $p_i$), however, has been
fairly well studied (e.g. Bayer and Diaconis $\cite{Bayer}$ or Diaconis,
McGrath, and Pitman $\cite{DiMcPi}$).

\section{Bounding the Time to Uniform}

	This section uses the concept of a strong uniform time to upper bound the
time for biased riffle shuffles to get close to the uniform
distribution. The bounds obtained are of the same order as lower bounds due to
Lalley \cite{Lalley}.

	Recall that the total variation
distance between two probability distributions $P_1$ and $P_2$ on a set $X$
is defined as:

\[ \| P_1-P_2 \| = \frac{1}{2} \sum_{x \in X} |P_1(x)-P_2(x)| \]

	Let $P^{*k}$ denote the $k$-fold convolution of $P$. Let $U$ be the
uniform distribution on $S_n$.

\begin{theorem} \label{upper}

\[ \| P_{n,a;\vec{p}}^{*k}-U \| \leq {n \choose 2} [p_1^2+\cdots+p_a^2]^k \]

\end{theorem}

\begin{proof}
	For each $k$, let $A^k$ be a random $n*k$ matrix formed by letting
each entry equal $i$ with probability $p_i$. Note that the random matrix
$A^k$ corresponds to a random permutation under the measure
$P^{*k}_{n,a;\vec{p}}$. To see this, recall Description 4 of biased riffle
shuffles (the inverse description). A single inverse $a$ shuffle
corresponds to a column of $A^k$ by letting the $i$th entry in the column
of $A^k$ equal the pile into which card $i$ is placed.

	Let $T$ be the first time that the rows of $A^k$ are
distinct. It is not hard to see that $T$ is a strong uniform time for
$P^{*k}_{n,a;\vec{p}}$ in the sense of Sections 4B-4D of Diaconis
$\cite{Diac}$. Namely, the permutation associated to the matrix $A^T$
is uniform. This is because, as in Proposition $\ref{convol}$, the
inverse of the $k$ fold convolution of $a$-shuffles may be viewed as
inverse sorting into $a^k$ piles, and at time $T$ each pile has at
most 1 card. Symmetry implies that these cards are in uniform random
order. It is proved on page 76 of Diaconis $\cite{Diac}$ that:

	\[ |P^{*k}_{n,a;\vec{p}}-U| \leq Prob(T > k) \]

	Let $V_{ij}$ be the event that rows $i$ and $j$ of $A^k$ are
the same. The probability that $V_{ij}$ occurs is
$[p_1^2+\cdots+p_a^2]^k$. The theorem follows since:

\begin{eqnarray*}
Prob(T>k) & = & Prob (\cup_{1 \leq i < j \leq n}) A_{ij}\\
& \leq & \sum_{1 \leq i < j \leq n} Prob(A_{ij})\\
& = & {n \choose 2} [p_1^2+\cdots+p_a^2]^k
\end{eqnarray*}

\end{proof}

{\bf Remarks}

\begin{enumerate}

\item Theorem $\ref{upper}$ shows that $k=2log_{\frac{1}{\sum_{i=1}^a p_i^2}}n$
steps suffice to get close to the uniform distribution (in the case
$a=2,p_1=p_2=\frac{1}{2}$ this is $2log_2n$).

\item Lalley \cite{Lalley} proved that there exists an open neighborhood of
$p_1=\frac{1}{2}$ such that for all $p_1$ in this neighborhood,
a $P_{n,2;p_1,p_2}$ shuffle takes at least

	\[ \frac{3+\theta}{4} log_{\frac{1}{p_1^2+p_2^2}}n \]

	steps to get close to the uniform distribution. Here $\theta=\theta_{p_1}$ is
the unique real number such that

	\[ p_1^{\theta}+p_2^{\theta} = (p_1^2+p_2^2)^2 \]

	Note that when $p_1=p_2=\frac{1}{2}$ this bound is $\frac{3}{2}log_2 n$, which
is of the same order as the $2 log_2 n$ bound of Theorem \ref{upper}, and agrees
exactly with the more refined analysis of Bayer and Diaconis \cite{Bayer} for the
GSR shuffles.

\end{enumerate}

\section{Gessel's Bijection and Cycle Structure} \label{bijec}

	This section begins by describing a bijection of Gessel
$\cite{Gessel}$. This requires some preliminary notation and
concepts. Recall that a permutation $\pi \in S_n$ is said to have a descent
at position $i$ if $\pi(i) > \pi(i+1)$. We adopt the convention that all
$\pi \in S_n$ have a descent at position $n$. Define a necklace on an
alphabet to be a sequence of cyclically arranged letters of the alphabet. A
necklace is said to be primitive if it is not equal to any of its
non-trivial cyclic shifts. For example, the necklace $(a\ a\ b\ b)$ is
primitive, but the necklace $(a\ b\ a\ b)$ is not.

	Given a word $w$ of length $n$ on an ordered alphabet, the 2-row
form of the standard permutation $st(w) \in S_n$ is defined as
follows. Write $w$ under $1 \cdots n$ and then write under each letter of
$w$ its lexicographic order in $w$, where if two letters of $w$ are the
same, the one to the left is considered smaller. For example (page 195 of
Gessel and Reutenauer $\cite{Gessel}$):

\[ \begin{array}{c c c c c c c c c c c c c c}
	& & 1 & 2 & 3 & 4 & 5 & 6 & 7 & 8 & 9 & 10 & 11 & 12\\
	w & = & b & b & a & a & b & c & c & c & b & c & b & b\\
	st(w) &= & 3 & 4 & 1 & 2 & 5 & 9 & 10 & 11 & 6 & 12 & 7 & 8
	\end{array} \]

	For a finite ordered alphabet $A$, Gessel and Reutenauer
$\cite{Gessel}$ give a bijection $U$ from the set of length $n$ words $w$
of onto the set of finite multisets of necklaces of total size $n$, such
that the cycle structure of $st(w)$ is equal to the cycle structure of
$U(w)$. To define $U(w)$, one replaces each number in the necklace of
$st(w)$ by the letter above it. In the example, the necklace of $st(w)$ is
$(1\ 3), (2\ 4), (5), (6\ 9), (7\ 11\ 8\ 12\ 10)$. This gives the following
multiset of necklaces on $A$:

\[ (a\ b)(a\ b)(b)(b\ c)(b\ c\ b\ c\ c) \]
 
	Theorem $\ref{Gesbij}$, one of the main results of this section,
will follow from this bijection.

\begin{theorem} \label{Gesbij} Fix $r_1,\cdots,r_a \geq 0$ such that
$\sum_{i=1}^a r_i=n$. The bijection $U$ defines by restriction a
cycle-structure preserving bijection $\bar{U}$ from elements of $S_n$ with
descent set contained in $\{r_1, r_1+r_2 , \cdots, r_1+\cdots+r_a=n\}$ to
multisets of primitive necklaces on the alphabet $\{1,\cdots,a\}$ formed
from a total of $r_i$ $i$'s.  \end{theorem}

\begin{proof}
	Restrict $U$ to the set of words with $r_i$ $i's$. It is clear
that an element $\pi$ of $S_n$ can arise as the standard permutation
of at most one word with $r_i$ $i's$. Also, the $\pi$ which arise are
precisely those $\pi$ such that the descent set of $\pi^{-1}$ is
contained in $\{r_1, r_1+r_2 , \cdots, r_1+\cdots+r_a=n\}$. This
proves the theorem.
\end{proof}

	Corollary $\ref{translate}$ will translate Theorem $\ref{Gesbij}$
into the language of generating functions. This uses some further
notation. Define the quantity $M(r_1,\cdots,r_a)$ as:

\[ M(r_1,\cdots,r_a) = \frac{1}{n} \sum_{d|n,r_1,\cdots,r_a} \mu(d)
\frac{\frac{n}{d}!}{\frac{r_1}{d}! \cdots \frac{r_a}{d}!} \]

	One easily proves by Moebius inversion (e.g. page 172 of Hall
$\cite{Hall}$) that $M(r_1,\cdots,r_a)$ is the number of primitive
circular words from an alphabet $\{1,\cdots,a\}$ in which the letter
$i$ appears $r_i$ times.

	Recall that we are using the convention that all permutations
in $S_n$ have a descent at position $n$. For $b_i,n_i \geq 0$, let
$\vec{b}=(b_1,\cdots,b_a)$ and $\vec{n}=(n_1,n_2,\cdots)$. Let
$A_{\vec{b},\vec{n}}$ be the number of permutations on $b_1 + \cdots +
b_a$ letters with descent set contained in $\{b_1,b_1+b_2,\cdots,
b_1+\cdots+b_a\}$ and $n_i$ $i$-cycles.

\begin{cor} \label{translate} For all $a \geq 1$,

\[ \sum_{\vec{b},\vec{n}} A_{\vec{b},\vec{n}} \prod_{i=1}^a z_i^{b_i}
\prod_j x_j^{n_j} = \prod_{i=1}^{\infty} \prod_{r_1,\cdots,r_a \geq 0
\atop r_1+\cdots+r_a=i} (\frac{1}{1-z_1^{r_1}\cdots
z_a^{r_a}x_i})^{M(r_1,\cdots,r_a)} \]
 
\end{cor}

\begin{proof}
	The coefficient of $\prod_{i=1}^a z_i^{b_i} \prod_j x_j^{n_j}$
on the left hand side is equal to $A_{\vec{b},\vec{n}}$, the number of
permutations on $b_1+\cdots+b_a$ letters with descent set contained in
$\{b_1,b_1+b_2,\cdots, b_1+\cdots+b_a\}$ and $n_j$ $j$-cycles. Theorem
$\ref{Gesbij}$ says that this is the number of multisets of necklaces
on the alphabet $\{1,\cdots,a\}$ with $b_i$ $i$'s and $n_j$
$j$-cycles. The corollary now follows from the interpretation of
$M(r_1,\cdots,r_a)$ as the number of primitive circular words of
length $n$ from an alphabet $\{1,\cdots,a\}$ in which the letter $i$
appears $r_i$ times.
\end{proof}

	Corollary $\ref{translate}$ will be used to study the cycle structure
of a permutation under the measure $P_{n,a,\vec{p}}$. Let $E_{n,a,\vec{p}}$
denote expectation with respect to the measure
$P_{n,a,\vec{p}}$, and $N_i$ denote the random variable on $S_n$ such that
$N_i(\pi)$ is the number of $i$-cycles of $\pi$. The case of Theorem
$\ref{cycle}$ with all $p_i=\frac{1}{a}$ is known from Diaconis, McGrath,
and Pitman $\cite{DiMcPi}$.

\begin{theorem} \label{cycle}
 
\[ \sum_{n=0}^{\infty} u^n E_{n,a,\vec{p}} \prod_{i=1}^N x_i^{N_i} =
\prod_{i=1}^{\infty} \prod_{r_1,\cdots,r_a \geq 0 \atop r_1+\cdots+r_a=i}
(\frac{1}{1-p_1^{r_1} \cdots p_a^{r_a} u^i x_i})^{M(r_1,\cdots,r_a)} \]

\end{theorem}

\begin{proof}
	Corollary $\ref{translate}$ and elementary manipulations imply
that:

\begin{eqnarray*}
\prod_{i=1}^{\infty} \prod_{r_1,\cdots,r_a \geq 0 \atop r_1+\cdots+r_a=i}
(\frac{1}{1-p_1^{r_1} \cdots p_a^{r_a} u^i x_i})^{M(r_1,\cdots,r_a)} & = & 
\sum_{n=0}^{\infty} u^n \sum_{b_1+\cdots+b_a=n \atop \vec{n}: \sum in_i=n} A_{\vec{b},\vec{n}} \prod_{i=1}^a p_i^{b_i} \prod_j x_j^{n_j}\\
& = & \sum_{n=0}^{\infty} u^n \sum_{b_1+\cdots+b_a=n \atop \vec{n}: \sum in_i=n} [{n
\choose b_1 \cdots b_a} \prod_{i=1}^a p_i^{b_i}]
[\frac{A_{\vec{b},\vec{n}}}{{n \choose b_1 \cdots b_a}}] \prod_j x_j^{n_j}
\end{eqnarray*}

	We give a probabilistic interpretation to:

\[ \sum_{n=0}^{\infty} u^n \sum_{b_1+\cdots+b_a=n \atop \vec{n}: \sum
in_i=n} [{n \choose b_1 \cdots b_a} \prod_{i=1}^a p_i^{b_i}]
[\frac{A_{\vec{b},\vec{n}}}{{n \choose b_1 \cdots b_a}}] \prod_j
x_j^{n_j} \]

	The first term in square brackets is the chance that a deck cut
according to the $mult(n,\vec{p})$ distribution is cut into packets of size
$b_1,\cdots,b_a$. To interpret the second term in square brackets, use the
fact from page 17 of Stanley $\cite{Stanley}$ that the total number of
permutations on $n=b_1+\cdots+b_a$ letters with descent set contained in
$\{b_1,b_1+b_2,\cdots,b_1+\cdots+b_a\}$ is the multinomial coefficient ${n
\choose b_1 \cdots b_a}$. Thus the second term is equal to the chance that
choosing uniformly among permutations on $n$ letters whose inverse has
descent set contained in $\{b_1,b_1+b_2,\cdots,b_1+\cdots b_a\}$ gives a
permutation with $n_i$ $i$-cycles. This proves the theorem.
\end{proof}

	As an example of an application of Theorem $\ref{cycle}$, one obtains
an expression for the expected number of fixed points after a $k$-fold
convolution of the measure $P_{n,a,\vec{p}}$.

\begin{cor} \label{fixpoint} The expected number of fixed points of a
permutation under the $k$-fold convolution of $P_{n,a,\vec{p}}$ is:

	\[ \sum_{j=1}^n [p_1^j+\cdots+p_a^j]^k \]

\end{cor}

\begin{proof}
	Recall from the introductory section that the $k$-fold
convolution of an $a$-shuffle with parameters $(p_1,\cdots,p_a)$ is
equivalent to an $a^k$ shuffle with parameters equal to the $a^k$
possible products $p_{s_1} \cdots p_{s_k}$ where each $s_i \in
\{1,\cdots,a\}$ and repetition is allowed. Thus it suffices to prove
the corollary in the case $k=1$.

	In the generating function of Theorem $\ref{cycle}$, one wants
to set $x_1=x$, $x_i=1$ for $i \geq 2$, then differentiate with
respect to $x$, set $x=1$, and finally take the coefficient of $y^n$.

	Setting $x_1=x$, $x_i=1$ for $i \geq 2$ in the generating
function of Theorem $\ref{cycle}$ gives:

\[ \frac{1}{1-y} \frac{1-p_1y}{1-p_1xy} \cdots
\frac{1-p_ay}{1-p_axy} \]

	because the $x_1=x$ term contributes $\frac{1}{\prod_{i=1}^a
(1-p_ixy)}$ and the $x_i=1$ for $i \geq 2$ term contributes
$\frac{\prod_{i=1}^a (1-p_iy)}{1-y}$. The corollary now follows by
easy algebra.
\end{proof}

{\bf Remarks}

\begin{enumerate}

\item In the case of $p_i=\frac{1}{a}$, Corollary $\ref{fixpoint}$ shows
that the expected number of fixed points after $k$ $a$-shuffles is:

\[ \sum_{j=1}^n \frac{1}{a^{(j-1)k}} \]

which is known from Diaconis, McGrath, and Pitman $\cite{DiMcPi}$. In fact
Holder's inequality gives:

\[ \frac{1}{a^{j-1}} \leq p_1^j+\cdots+p_a^j \]

so that the expected number of fixed points is smallest for unbiased
riffle shuffles.

\item It turns out that for $\frac{1}{(p_1^2+\cdots+p_a^2)^k} \gg 1$, the
number of fixed points is close to its Poisson(1) limit. In fact fixed
points (and more generally other functionals of cycle structure) approach
their limit distribution more quickly than $P_{n,a,\vec{p}}$ approaches its
uniform limit.

\end{enumerate}

\section{Enumerative Applications of Gessel's Bijection}

	This section considers some enumerative applications of Theorem
$\ref{Gesbij}$. To begin, formulas will be found for
the chance that an
$n$-cycle in $S_n$ has a given descent set $J$. Recall that all
permutations in $S_n$ are considered to have a descent at position $n$. We also
use the notation that if
$J=\{j_1<j_2< \cdots j_d=n \}$ and $j_0=0$, then $C(J)$, the composition of
the descent set $J$, is equal to $(j_1-j_0,\cdots,j_d-j_{d-1})$.

	Stanley $\cite{Stanley}$ gives two formulas for the number of
permutations with descent set $J$. These will both turn out to have analogs
for the case of $n$-cycles.

\begin{prop} \label{Stan1} (Page 69 of Stanley $\cite{Stanley}$) The number of
elements of $S_n$ with descent set $J$ is:
                               
	\[ \sum_{K \subseteq J} (-1)^{|J|-|K|} {n \choose C(K)} \]

\end{prop}

	This carries over to $n$-cycles as follows, where
$M(r_1,\cdots,r_a)$ is defined as in Section $\ref{bijec}$. 

\begin{cor} \label{ncyc1} The number of $n$-cycles with descent set $J$ is:

\[ \sum_{K \subseteq J} (-1)^{|J|-|K|} M(C(K)) \]
\end{cor}

\begin{proof}
	By Moebius inversion on the power set of $\{1,\cdots,n\}$, it
suffices to show that the number of $n$ cycles with descent set
contained in $K$ is $M(C(K))$. This follows from Theorem
$\ref{Gesbij}$.
\end{proof}

	There is also a determinantal formula for the number of
permutations with descent set $J$. Suppose that the elements of $J$ are $1
\leq j_1 \leq j_2 ... \leq j_k \leq n-1$. Define $j_0=0$ and $j_{k+1}=n$.

\begin{prop} \label{Stan2} (Page 69 of Stanley $\cite{Stanley}$) The number of
elements of $S_n$ with descent set $J$ is the determinant of a $k+1$ by
$k+1$ matrix, where $(l,m) \in [0,k] \times [0,k]$:

	\[ det {n-j_l \choose j_{m+1}-j_l} \]

\end{prop}                               

	This can be generalized to $n$-cycles. Given $J$, a subset of
$\{1,...,n-1\}$, let $J^d$ be the subset of $J$ consisting of all numbers
divisible by $d$. If $J$ is non-empty, label these elements $1 \leq j^d_1
\leq j^d_2 \cdots \leq j^d_{|J^d|} \leq n-1$. Define $j^d_0=0$ and
$j^d_{|J^d|+1}=n$.

\begin{theorem} \label{ncycdet} The number of $n$-cycles
with descent set $J$ is:

\[ \frac{1}{n} \sum_{d|n} \mu (d) (-1)^{|J|-|J^d|} det
{\frac{n}{d}-\frac{j^d_l}{d} \choose
\frac{j^d_{m+1}}{d}-\frac{j^d_l}{d}} \]

\end{theorem}

\begin{proof}
	From Theorem $\ref{ncyc1}$, the number of $n$-cycles with
descent set $J$ is:

\begin{eqnarray*}
\sum_{K \subseteq J} (-1)^{|J|-|K|} M(C(K)) & = & \frac{1}{n} \sum_{K
\subseteq J} (-1)^{|J|-|K|} \sum_{d: K \subseteq J^d} \mu(d) {\frac{n}{d} \choose
C(\frac{K}{d})}\\
& = & \frac{1}{n} \sum_{d|n} \mu(d) \sum_{K \subseteq J^d} (-1)^{|J|-|K|}
{\frac{n}{d} \choose C(\frac{K}{d})}\\ & = & \frac{1}{n} \sum_{d|n} \mu(d)
(-1)^{|J|-|J^d|} \sum_{K \subseteq J^d} (-1)^{|J^d|-|K|} {\frac{n}{d} \choose
C(\frac{K}{d})}
\end{eqnarray*}
	
	Proposition $\ref{Stan1}$ shows that $\sum_{K \subseteq J^d}
(-1)^{|J^d|-|K|} {\frac{n}{d} \choose C(\frac{K}{d})}$ is the
number of permutations on $\frac{n}{d}$ symbols with descent set
$\frac{J^d}{d}$. The theorem then follows from Proposition
$\ref{Stan2}$.
\end{proof}

	The enumeration of matrices with fixed row and column sums is
related to some problems in statistics (see for instance the work of
Diaconis and Sturmfels $\cite{DiaSturm}$). Proposition $\ref{invol}$
relates the theory of such matrices to the theory of descents in
involutions.

\begin{prop} \label{invol} The number of involutions in $S_n$ with descent
set contained in $K=\{k_1,...,k_r=n\}$ is equal to the number of symmetric $r*r$
matrices with non-negative integer entries and with $i$th row sum
$k_i-k_{i-1}$, where by convention $k_0=0$.  \end{prop}

\begin{proof}
	Theorem $\ref{Gesbij}$ shows that it suffices to count the number
of multisets of primitive necklaces on an alphabet of $k_i-k_{i-1}$ $i$'s,
where each necklace has length 1 or 2. Note that a primitive necklace of
length 2 consists of a pair of distinct elements. So for $i \neq j$, let
$X_{ij}$ be the number of pairs of letter i with letter j, and
let $X_{ii}$ be the number of singleton $i$'s. The matrix $(X_{ij})$ has
all the desired properties.
\end{proof}

\section{Inversion and Descent Structure After a Shuffle}

	It is natural to study the inversion and descent structure of a
permutation obtained after a biased riffle shuffle.	Recall that $\pi$ is said to
invert the pair $(i,j)$ with $i<j$ if
$\pi(i)>\pi(j)$. The number of inversions of $\pi$ is the number of pairs
which $\pi$ inverts and will be denoted $Inv(\pi)$. It is easy to see that
$Inv(\pi)=Inv(\pi^{-1})$ and that $Inv(\pi)$ is the length of $\pi$ in
terms of the generators $\{(i,i+1):1 \leq i \leq n-1 \}$. Theorem
$\ref{invgen}$ will give a $q$-exponential generating function for $Inv$
after a biased riffle shuffle. This uses the notation:

\[ [n]! = \prod_{i=0}^{n-1} (1+q+\cdots+q^i) \]

\[ \qb{n}{k} = \frac{[n]!}{[k]![n-k]!} \]

	As usual, $E_{n,a,\vec{p}}$ denotes expectation with respect to the
measure $P_{n,a,\vec{p}}$. As will be explained in the course of the proof,
the second equality in Theorem $\ref{invgen}$ is purely formal in the sense
that it only holds if $|q|<1$, and thus only the first equality should be used for
the purpose of computing moments.

\begin{theorem} \label{invgen}

\begin{eqnarray*}
\sum_{n=0}^{\infty} \frac{u^n}{[n]!}
E_{n,a,\vec{p}} \  q^{Inv} & = & \prod_{i=1}^a [\sum_{j=0}^{\infty}
\frac{(up_i)^j}{[j]!}] \\
& = & \prod_{i=1}^a \prod_{j=0}^{\infty} \frac{1}{1-up_i(1-q)q^j}
\end{eqnarray*}

\end{theorem} 

\begin{proof}
	The following identity is clear from elementary manipulations and
the definition of $q$-multinomial coefficients:

\[ \sum_{n=0}^{\infty} \sum_{b_i \geq 0 \atop b_1+\cdots+b_a=n}
\frac{p_1^{b_1}\cdots p_a^{b_a} \qb{n}{b_1 \cdots b_a} u^n}{[n]!} =
\prod_{i=1}^a [\sum_{j=0}^{\infty} \frac{(up_i)^j}{[j]!}] \]

	The left-hand side can be rewritten as:

\[ \sum_{n=0}^{\infty} \frac{u^n}{[n]!} \sum_{b_i \geq 0 \atop
b_1+\cdots+b_a=n} [ {n \choose b_1 \cdots b_a} \prod_{i=1}^a p_i^{b_i}]
\frac{\qb{n}{b_1 \cdots b_a}}{ {n \choose b_1 \cdots b_a}} \]

	Since $Inv(\pi)$ is equal to $Inv(\pi^{-1})$, it is sufficient to
analyze the number of inversions in the inverse of a permutation chosen
from the measure $P_{n,a;\vec{p}}$. Recalling the first description of
biased riffle shuffling in Section $\ref{intro}$, note that the term in
brackets corresponds to picking the packet sizes $b_1,\cdots,b_a$ according
to the $mult(a;\vec{p})$ law. From pages 22 and 70 of Stanley
$\cite{Stanley}$, it is known that $\qb{n}{b_1 \cdots b_a}$ is the sum of
$q^{Inv(\pi)}$ over all $\pi$ in $S_n$ with descent set contained in
$\{b_1,b_1+b_2,\cdots,b_1+\cdots+b_a=n\}$ and that $ {n \choose b_1 \cdots
b_a}$ is the number of permutations with descent set contained in
$\{b_1,b_1+b_2,\cdots,b_1+\cdots+b_a=n\}$. These observations prove the
first equality of the theorem.

	The second equality follows from a famous identity of Euler, which
is true if $|x|,|q|<1$:

\[ \prod_{j=0}^{\infty} \frac{1}{1-xq^n} = \sum_{j=0}^{\infty}
\frac{x^j}{(1-q) \cdots (1-q^j)} \]
\end{proof}

	Theorem $\ref{invgen}$ can be used to compute the expected number
of inversions after a $k$-fold convolution of a $P_{n,a;\vec{p}}$
shuffle. However, we prefer the following direct probabilistic argument.

\begin{theorem} \label{expecinversion}

	The expected number of inversions under the $k$-fold convolution
of $P_{n,a;\vec{p}}$ is:

\[ \frac{{n \choose 2}}{2} [1-(p_1^2+\cdots+p_a^2)^k] \]

\end{theorem}

\begin{proof}
	For $1 \leq i<j \leq n$, define a random variable $X_{i,j}$ as follows. In the
inverse model of card shuffling, let $X_{i,j}=1$ if card $i$ goes to a pile to
the right of card $j$, and let $X_{i,j}=0$ otherwise. It is easy to see that if
$\pi$ is the permutation obtained after the shuffle, then $\pi(i)>\pi(j)$ exactly
when $X_{i,j}=1$. Thus,

\[ Inv = \sum_{1 \leq i<j \leq n} X_{i,j}. \]

	It is clear that each $X_{i,j}$ has expected value $\frac{1-(p_1^2+\cdots+
p_a^2)^k}{2}$, because this is one half the chance that cards $i$ and $j$ fall in
different piles. The theorem now follows by linearity of expectation.
\end{proof}

{\bf Remarks}

\begin{enumerate}

\item Note that a uniformly chosen element of $S_n$ has on average $\frac{{n
\choose 2}}{2}$ inversions. In fact the distribution for inversions on $S_n$ is
the sum $X_1+\cdots+X_n$ where the $X_i$ are independent and uniform on $[0,i-1]$.

\item By Holder's inequality, the expected number of inversions is maximum
for $k$ unbiased $a$ shuffles (which is the same as an $a^k$ shuffle), and
in this case is $\frac{{n \choose 2}}{2} [1-\frac{1}{a^k}]$. For instance,
a $1$ shuffle of a sorted deck gives no inversions, and a 2 shuffle of a
sorted deck gives a permutation which has on average one half as many
inversions as a random permutation.

\item	It would be interesting to use Theorem $\ref{invgen}$ to study the
asymptotics of inversions after a biased riffle shuffle. Even for the case
$a=2,p_1=p_2=\frac{1}{2}$, it is not known if the $n \rightarrow \infty$
limit distribution is normal.

\item The same technique used in Theorem \ref{expecinversion} can be used to
study the distribution of $Des(\pi)$, the number of descents of a permutation
$\pi$ after a biased riffle shuffle. For example, using the convention that all
elements of $S_n$ have a descent at position $n$, the expected number of descents
would be

	\[ 1 + \frac{n-1}{2} [1-(p_1^2+\cdots+p_a^2)^k] \]

	It is perhaps surprising that these moments can be computed so easily. One
reason to be surprised is that in the case of unbiased shuffles, Bayer and
Diaconis \cite{Bayer} showed that $Des(\pi^{-1})$ is a sufficient statistic for
the random walk. Nevertheless, computing the moments of $Des(\pi^{-1})$ is
more difficult than computing the moments of $Des(\pi)$, as a glance at the work
of Mann \cite{Mann} will make clear.

\end{enumerate}

\section{Acknowledgements} The author is grateful to Persi Diaconis for numerous
conversations about the theory of descents and shuffling cards. Ira Gessel was
also kind enough to share his thoughts on descents. This research was done under
the support of the National Defense Science and Engineering Graduate Fellowship
(grant no. DAAH04-93-G-0270) and the Alfred P. Sloan Foundation Dissertation
Fellowship.

\end{document}